\documentclass[journal,a4paper]{IEEEtran}%

\usepackage{makeidx}
\usepackage{graphicx}
\usepackage{mathptmx}
\usepackage{amsmath}
\usepackage{latexsym}
\usepackage{epsfig}
\usepackage{color}
\usepackage{algorithm}
\usepackage{algorithmicx}

\begin{document}
\title{$N$-$k$-$\epsilon$ Survivable Power System Design}

\author{Richard Li-Yang Chen, Amy Cohn, Neng Fan, Ali Pinar
\thanks{R. Chen and A. Pinar are with Quantitative Modeling and Analysis, Sandia National Laboratories, Livermore, California 94551, USA. email: \{rlchen,apinar\}@sandia.gov.}
\thanks{A. Cohn is with Industrial and Operations Engineering, University of Michigan, Ann Arbor, Michigan 48109, USA. email: amycohn@umich.edu.}
\thanks{N. Fan is with Discrete Math \& Complex Systems, Sandia National Laboratories, Albuquerque, New Mexico 87185, USA. email: nnfan@sandia.gov.}}

\maketitle

\begin{abstract}
We consider the problem of designing (or augmenting) an electric power system such that it satisfies the $N$-$k$-$\epsilon$ survivability criterion while minimizing total cost. The survivability criterion requires that at least $(1-\epsilon)$  fraction of the total demand can still be met even if any $k$ (or fewer) of the system components fail. We formulate this problem, taking into account both transmission and generation expansion planning, as a mixed-integer program.  Two algorithms are designed and tested on modified instances from the IEEE-30-Bus and IEEE-57-Bus systems.
\end{abstract}

\begin{IEEEkeywords}
Long-term grid planning, contingency requirements, decomposition, separation oracle, implicit optimization.
\end{IEEEkeywords}

\maketitle

\section{Introduction}\label{sec1}
\IEEEPARstart{A}{ccording} to the Transmission Planning Standard (TPL-001-1, \cite{NERC}), defined by the North American Electric Reliability Corporation (NERC), power systems are required to perform necessary adjustments under normal and contingency conditions to ensure system reliability. If  only a single element is lost ($N$-$1$ contingency), the system must be stable and all thermal and voltage limits must remain within applicable rating. The loss-of-load is not allowed for $N$-$1$ contingency. In the case of multiple simultaneous  failures ($N$-$k$ contingency), the system still has to meet the stable, thermal and voltage limits, but planned or controlled loss-of-load is allowed, to a limited degree.

Recently, optimization methods have been applied to $N$-$k$ contingency analysis for a large variety of vulnerabilities within power systems. For example, line-vulnerability studies can be found in \cite{Pinar2010,Bienstock2010,Arroyo2010}, where optimization methods are used to find small groups of lines whose failure can cause severe blackout or large loss-of-load. $N$-$k$ contingencies are also considered in optimal power flow models \cite{SalWB04,SalWB09,Fan2011} and unit commitment problems \cite{Street2011}. The methods used in \cite{Arroyo2010,Fan2011,Street2011} are all based on a bilevel programming approach, which is the main method used for network inhibition/interdiction problems.

For power system expansion problem with the added consideration of contingencies, \cite{Romero2011} proposed a multilevel mixed integer programming model and solved it by tabu search.  They use this model to analyze the interaction between a power system defender and a terrorist, who seeks to disrupt system operations. For the transmission expansion problem, references \cite{Carrion2007,Choi2007,Moulin2010} considered the contingency criteria by stochastic programming and integer programming approaches. The generation expansion problem has recently been studied in \cite{Jin2011}, which pointed to many recent advances. Transmission and generation expansion planning problems have also been studied in the context of renewable energy integration, see, e.g., \cite{Bent2011}.

In this paper, we consider the transmission and generation expansion planning (TGEP) problem of designing (or augmenting) an electric power system of minimum total cost that satisfies the $N$-$k$-$\epsilon$ survivability criterion. The survivability criterion requires that a feasible power flow must still exist, satisfying at least $(1-\epsilon)$ fraction of the total demand, even after failures of up to $k$ elements of the power system. Considering the standards of NERC, for no-contingency state and contingency states with $k=1$, no loss-of-load is allowed; for contingency states with $k\ge 2$, a small fraction of total load demand can be shed.

We formulate a mixed-integer nonlinear program (MINLP) to model TGEP along with multiple states representing all the possible contingency scenarios and the corresponding flow variables to ensure that $(1-\epsilon)$ fraction of the demand can be met. The combinatorial number of contingency scenarios imposes a substantial computational burden. To overcome this challenge, we propose two cutting plane algorithms, one based on a Benders decomposition method to check the load satisfaction of each contingency scenario and another based on a custom cutting plane algorithm, which solves a bilevel separation problem to determine the worst-case loss-of-load under any contingency with up to $k$ failures. To test our models and algorithms, numerical experiments are performed on the IEEE-30-Bus system and the IEEE-57-Bus system.

The rest of this paper is organized as follows: In Section \ref{sec2}, the TGEP problem considering the full set of contingency scenarios is formulated as a MINLP; Section \ref{sec3} presents two methods to solve this large-scale MINLP; In Section \ref{sec4}, numerical experiments are performed on two IEEE test systems; Section \ref{sec5} concludes the paper.

\section{Models}\label{sec2}
\subsection{Nomenclature}
\noindent Sets and indices
\begin{tabbing}
\hspace{6ex}\=\hspace{6ex}\=\hspace{4ex}\=\hspace{4ex}\=\hspace{4ex}\=\hspace{4ex}\=\hspace{4ex}\=\hspace{4ex}\kill
$I$ \>      Set of buses (indexed by $i,j$).\\
$S(k)$\>  Set of \emph{ALL} contingency states with $k$ or fewer failures.\\
$G$ \>  Set of generating units.\\
$G_i$\> Set of generating units at bus $i$.\\
$E$ \>  Set of transmission elements.\\
$E_{.i}$ \>  Set of transmission elements oriented into bus $i$.\\
$E_{i.}^s$ \>  Set of transmission elements oriented out of bus $i$.\\
$i_e,j_e$ \>    Tail/head (bus no.) of transmission element $e=(i_e,j_e)$.\\
\end{tabbing}

\noindent Parameters
\begin{tabbing}
\hspace{6ex}\=\hspace{6ex}\=\hspace{4ex}\=\hspace{4ex}\=\hspace{4ex}\=\hspace{4ex}\=\hspace{4ex}\=\hspace{4ex}\kill
$C_e$ \>     Investment cost of transmission element $e$.\\
$C_g$\>      Investment cost of generating unit $g$.\\
$C^p_g$ \>   Marginal production cost of generating unit $g$.\\
$\overline{P}_g$ \>   Maximum capacity of unit $g$.\\
$B_e$ \>     Electrical susceptance of transmission element $e$.\\
$F_{e}$ \>   Capacity of transmission element $e$.\\
$D_i$ \>     Electricity load demand at bus $i$.\\
$\sigma$ \>  Weighting factor to make investment cost and operating\\
\> cost comparable.\\
$\epsilon$\>    Fraction of load demand that can be shed.\\
$k$ \>      Contingency budget indicating the maximum number\\
\> of failed elements.\\
$\hat d_g^s$ \>        Binary parameter that takes value 1 if generating unit $g$ \\
\>                is part of the contingency state $s$ and 0 otherwise.\\
$\hat d_e^s$ \>        Binary parameter that takes value 1 if transmission  \\
\>                element $e$ is part of the contingency state $s$ and 0\\
\> otherwise.
\end{tabbing}

\noindent Decision Variables
\begin{tabbing}
\hspace{6ex}\=\hspace{6ex}\=\hspace{4ex}\=\hspace{4ex}\=\hspace{4ex}\=\hspace{4ex}\=\hspace{4ex}\=\hspace{4ex}\kill
$x_g$ \>  Binary generation expansion variable that takes value\\
\>        1 if  generating unit $g$ is added and 0 otherwise.\\
$x_e$ \>  Binary transmission expansion variable that takes value \\
\>        1 if  transmission element $e$ is added and 0 otherwise.\\
$q^{s}_i$ \>  Loss of load at bus $i$ for state $s$.\\
$p^{s}_g$ \>  Power output of generating unit $g$ for state $s$.\\
$f^{s}_e$ \>  Power flow  for transmission element $e$ for state $s$.\\
$\theta_i^{s}$ \> Phase angle of bus $i$ for state $s$.\\
$d_g^s$ \>        Binary variable that takes value 1 if generating unit $g$ \\
\>                is part of the contingency state $s$ and 0 otherwise.\\
$d_e^s$ \>        Binary variable that takes value 1 if transmission  \\
\>                element $e$ is part of the contingency state $s$ and 0\\
\> otherwise.
\end{tabbing}

For a contingency state $s\in S(k)$, $\hat{d}_g^s=1$ and $\hat{d}_e^s=1$ denote generating unit $g$ and transmission element $e$ fails in state $s$, respectively. Conversely, $\hat{d}_g^s=0$ and $\hat{d}_e^s=0$ denote that both these two elements are available. Therefore, in the no-contingency state $(s=0)$, $\hat{d}_g^0=0$ and $\hat{d}_e^0=0$ for all $g\in G$ and $e\in E$.

\subsection{TGEP Model}
In the following model, we extend the standard TGEP problem to include contingency constraints. Without loss of generality, we treat all power system elements as candidates; for an existing element, the investment cost $(C_e, C_g)$ can simply be set to $0$ and the corresponding investment decision $(x_e, x_g)$ is fixed to be $1$.

Once the planning decision is made, each newly planned element is available in all contingencies scenario $s\in S(k)$, unless the it is part of a given contingency. Additionally, in the no-contingency state $(s=0)$ no loss-of-load is allowed. For contingency state $s\in S(k) (s>0)$, the total loss-of-load is limited by the threshold $\epsilon$. The MINLP model for TGEP is formulated as follows,
\begin{subequations}\label{ef}
\begin{align}
\min_{x,f,p,q,\theta}  \ &\sum_{e \in E} C_e x_e + \sum_{g\in G} C_gx_g + \sigma  \sum_{g\in G}  C_g^p p_g^0 \label{ef_obj} \\
\textrm{s.t.} \quad
	&\sum_{g\in  G_i^s} p_g^{s}  + \sum_{e \in E_{.i}^s} f_{e}^{s} -\sum_{e \in E_{i.}^s} f_{e}^{s} + q_i^s = D_i, \quad  \forall i, s\label{ef_node_bal} \\
    &B_{e}\big(\theta_{i_e}^{s}- \theta_{j_e}^{s}\big) x_e (1-\hat d_e^s) - f_{e}^{s}=0, \quad \forall e, s \label{ef_flow_cons}\\
	&-F_e x_e(1-\hat d_e^s)\leq f_{e}^{s} \leq  F_{e}x_e(1-\hat d_e^s), \quad \forall e, s  \label{ef_flow_cap} \\
	&0 \leq p_g^{s} \leq \overline{P}_g x_g(1-\hat d_g^s), \quad \forall g, s  \label{ef_p_bounds}\\
    &0 \leq q^{s}_i \leq D_i, \quad \forall i, \forall s \in S(k)\setminus 0 \label{ef_q_lb}\\
    &\sum_{i \in I}  q_i^s \le \epsilon \sum_{i \in I}D_i,  \quad \forall s\in S(k)\setminus 0 \label{ef_lolc} \\
	&q^{0}_i = 0, \quad \forall i \label{ef_q_lb0}\\
	&x_g \in \{0,1\}, \quad \forall g \label{ef_x}\\
	&x_e \in \{0,1\}, \quad \forall e \label{ef_x}
\end{align}
\end{subequations}
In all subsequent formulations, unless otherwise specified, the indices $i,g,e$ and $s$ are elements of their corresponding sets, i.e., $i\in I, g\in G, e\in E$ and $s\in S(k)$.

The objective \eqref{ef_obj} is to minimize the total transmission and generation investment cost plus  the normalized operating cost in the no-contingency state $(s=0)$.

Constraints \eqref{ef_node_bal} are flow balance requirements for each bus and contingency pair. For any transmission element that is operational, Kirchhoff's voltage law must be enforced by \eqref{ef_flow_cons}. Power flow on transmission element $e$ is governed by thermal capacity constraints \eqref{ef_flow_cap}. For each contingency state, the power output of a generating unit must satisfy the upper bound given by \eqref{ef_p_bounds}.

The set of all contingency states, $S(k)$ includes the no-contingency state $s=0$.  In the $s=0$ state, all power system elements are available, and no loss-of-load is allowed, as limited by \eqref{ef_q_lb0}.

For contingency state $s>0$, constraints \eqref{ef_q_lb} and \eqref{ef_lolc} define the loss-of-load at each bus and across all possible contingencies, respectively. For contingency state with one failed element, no load shedding is allowed (i.e., $\epsilon=0$).  For states with two or more failed elements, $\epsilon > 0$ and the bound $ \epsilon \sum_{i \in I}D_i$ limits the total loss-of-load in the system. Therefore, for every contingency, at least $(1-\epsilon)\sum_{i \in I}D_i$ of  demand must be satisfied.

Observe that constraints \eqref{ef_node_bal}-\eqref{ef_q_lb} are specific to a particular contingency state, that is, for a given contingency state $s$, the transmission and generation elements in the contingency have zero capacity.  In the no-contingency state $s=0$, all invested transmission elements and generating units are available for the power flow problem.

\section{Solution Approaches}\label{sec3}
Replacing constraints \eqref{ef_flow_cons} by
\begin{align}
&B_{e}\big(\theta_{i_e}^{s} - \theta_{j_e}^{s}\big) - f_{e}^{s} + M_e(1-x_e +\hat d_e^s) \geq  0, \label{ef_flow_cons_lb}\\
&B_{e}\big(\theta_{i_e}^{s} - \theta_{j_e}^{s}\big) - f_{e}^{s} - M_e(1 -x_e + \hat d_e^s) \leq  0, \label{ef_flow_cons_ub}\\
& \hskip 4cm \forall e \in  E^s, s \in S(k), \nonumber
\end{align}
where $M_e$ is sufficiently large constant, formulation \eqref{ef} becomes a large-scale mixed integer linear program (MILP), which we refer as the \emph{extensive form} (EF).  EF  has an extremely large number of variables and constraints because it grows with the number of contingency states, which increases exponentially with $N$ and $k$.  For large power systems and/or  a contingency $k$ greater than one, EF  can quickly become computationally intractable.  In the following sections, we modify this formulation and present cutting plane algorithms for solving the reformulated problem.

\subsection{Benders Decomposition}
We begin by presenting an alternative formulation with only $|G|+|E|$ binary variables but possibly an extremely large number of constraints.  We use linear programming duality to generate valid inequalities for the projection of the natural formulation onto the space of the $x$ variables.  In essence, we use a variant of Benders Decomposition in which we generate valid inequalities corresponding to ``feasibility'' cuts.

For contingency state $s>0$, given a capacity expansion vector $\hat x$ and a contingency state vector $\hat{d}^s$, we solve the following linear program, denoted as the \emph{primal subproblem} PSP$(\hat x,\hat d^s)$, to determine an optimal power flow (OPF) that minimizes the loss-of-load.
\begin{subequations}\label{psp}
\begin{align}
   z(\hat x,\hat d^s) = &\min_{f,p,q,\theta} \ \sum_{i \in I} q_i^s \label{psp_obj} \\
	 \textrm{s.t. }     
	(\alpha_i^s)       \quad & \sum_{g\in  G_i^s} p_g^s+ \sum_{e \in E_{.i}^s} f_{e}^{s}-\sum_{e \in E_{i.}^s} f_{e}^{s}  +  q_i^s = D_i,  \  \forall i  \label{psp_flow_consv} \\
	(\hat \beta_e^s)   \quad &-B_{e}\big(\theta_{i_e}^{s} - \theta_{j_e}^{s}\big) + f_{e}^{s}  \leq M_e(1-\hat x_e + \hat d_e^s), \ \forall e   \\
	(\check \beta_e^s) \quad &B_{e}\big(\theta_{i_e}^{s} - \theta_{j_e}^{s}\big) - f_{e}^{s} \leq M_e(1-\hat x_e + \hat d_e^s), \ \forall e   \\
	(\delta_e^s)       \quad &f_{e}^{s} \leq  F_{e} \hat x_e(1-\hat d_e^s), \ \forall e     \label{psp_flow_ub} \\
	(\eta_e^s)         \quad &-f_{e}^{s} \leq  F_{e} \hat x_e(1-\hat d_e^s), \ \forall e     \label{psp_flow_lb} \\
	(\zeta_g^s)        \quad & 0  \leq p_g^{s} \leq \overline{P}_g \hat x_g(1-\hat d_g^s), \ \forall g  \label{psp_gen_ub} \\
	(\lambda_i^s)      \quad & 0  \leq q_i^{s} \leq D_i,\ \forall i   \label{psp_lol_ub}
\end{align}
\end{subequations}

The objective \eqref{psp_obj} is to minimize loss-of-load by adjusting the flow, phase angles and power generation, given the prescribed capacity expansion decision $\hat x$ and contingency $\hat d^s$.  Clearly, if $z(\hat x,\hat d^s)>\epsilon \sum_{i \in I} D_i$, there does not exist a feasible power flow satisfying at least $(1-\epsilon)$ of total demand, and if $z(\hat x,\hat d^s)\le \epsilon \sum_{i \in I} D_i$,  a feasible power flow exists that can satisfy \emph{at least} $(1-\epsilon)$ of total demand.

The variables indicated in parenthesis on the left-hand-side of the constraints in \eqref{psp} denote the corresponding dual variables. In turn, we can formulate the dual of this problem, DSP$(\hat x,\hat d^s)$ as follows,
\begin{align*}
\max_{\alpha, \hat \beta, \check \beta, \delta, \eta, \zeta, \lambda} \ \sum_{i \in I} D_i (\alpha_i^s + \lambda_i^s) + \sum_{e \in E} M_e(1-\hat x_e + \hat d_e^s)(\hat \beta_e^s + \check \beta_e^s) \\
+ \sum_{e \in E} F_e \hat x_e(1-\hat d_e^s)(\delta_e^s + \eta_e^s) + \sum_{g \in G} \overline{P}_g  \hat x_g(1-\hat d_g^s)\zeta_g^s,  \nonumber 
\end{align*}
subject to constraints corresponding to primal variables $f, p, q, \theta$. Since PSP$(\hat x,\hat d^s)$ has a finite optimal solution (in the worst case, all load will be shed), DSP$(\hat x,\hat d^s)$ also has a finite optimal solution and by strong duality, the optimal solutions coincide. Therefore, DSP$(\hat x,\hat d^s)$ has a finite optimal solution, and in fact, an optimal extreme point. Thus we can reformulate PSP$(\hat x,\hat d^s)$ as follows,
\begin{align}\label{dsp2}
&\max_{\ell=1,\cdots,L^s}  \sum_{i \in I} D_i (\alpha_i^{s\ell} + \lambda_i^{s\ell}) + \sum_{e \in E} M_e(1-\hat x_e + \hat d_e^s)(\hat \beta_e^{s\ell} + \check \beta_e^{s\ell}) \nonumber  \\
&\quad + \sum_{e \in E} F_e x_e(1-\hat d_e^s)(\delta_e^{s\ell} + \eta_e^{s\ell}) + \sum_{g \in G} \overline{P}_g  x_g(1-\hat d_g^s)\zeta_g^{s\ell},
\end{align}
where $L^s$ is the set of extreme points corresponding to the polyhedron characterized by dual constraints based on \eqref{psp} for primal variables $f, p, q, \theta$.

Observing that $z(\hat x,\hat d^s) \leq \epsilon \sum_{i \in I} D_i$ should be satisfied for all $s \in S(k)\setminus 0$, the contingency feasibility conditions can be defined as follows,
\begin{align}\label{dsp3}
& \sum_{i \in I} D_i (\alpha_i^{s\ell} + \lambda_i^{s\ell}) + \sum_{e \in E} M_e(1-\hat x_e + \hat d_e^s)(\hat \beta_e^{s\ell} + \check \beta_e^{s\ell})  \nonumber\\
&+ \sum_{e \in E} F_e \hat x_e(1-\hat d_e^s)(\delta_e^{s\ell} + \eta_e^{s\ell}) + \sum_{g \in G} \overline{P}_g  \hat x_g(1-\hat d_g^s)\zeta_g^{s\ell} \\
 &\hskip 3.7cm \leq \epsilon \sum_{i \in I} D_i, \  \forall \ell=1,\cdots,L^s \nonumber
\end{align}

Explicitly satisfying demand for the no-contingency state and using  \eqref{dsp3}  to satisfy the $(1-\epsilon)$ criterion for all contingency states with $k$ or fewer  failures, we can reformulate \eqref{ef} as:

\begin{subequations}\label{rmp}
\begin{align}
\min_{x,f,p,q,\theta}  &\  \sum_{e \in E} C_e x_e + \sum_{g\in G} C_gx_g + \sigma  \sum_{g\in G}  C_g^p p_g^0  \\
\textrm{s.t.} \ & \sum_{i \in I} D_i (\alpha_i^{s\ell} + \lambda_i^{s\ell}) + \sum_{e \in E} M_e(1- x_e + \hat d_e^s)(\hat \beta_e^{s\ell} + \check \beta_e^{s\ell})  \nonumber\\
	&+ \sum_{e \in E} F_e x_e(1-\hat d_e^s)(\delta_e^{s\ell} + \eta_e^{s\ell}) + \sum_{g \in G} \overline{P}_g  x_g(1-\hat d_g^s)\zeta_g^{s\ell} \nonumber \\
	&  \leq  \epsilon \sum_{i \in I} D_i, \ \forall \ell=1,\cdots,L^s, s \in S(k)\setminus 0 \label{rmp_f_cut} \\
	&\sum_{g\in  G_i} p_g^0+ \sum_{e \in E_{.i}} f_{e}^{0} - \sum_{e \in E_{i.}} f_{e}^{0}  =  D_i, \ \forall i \in  I  \\
    &B_{e}\big(\theta_{i_e}^{0} - \theta_{j_e}^{0}\big) - f_{e}^{0} + M_e(1-x_e+\hat d_e^0) \geq  0, \ \forall e\in E \\
    &B_{e}\big(\theta_{i_e}^{0} - \theta_{j_e}^{0}\big) - f_{e}^{0} - M_e(1-x_e+ \hat d_e^0) \leq  0, \ \forall e\in E \\
	&-F_e(1-\hat d_e^0) \leq f_{e}^{0} \leq  F_{e}(1-\hat d_e^0), \ \forall e \in  E   \\
	& 0 \leq p_g^{0} \leq \overline{P}_g x_g(1-\hat d_g^0), \ \forall g \in G  \\
	&x_g \in \{0,1\}, \ \forall g\in G \label{pf_x_g}\\
	&x_e \in \{0,1\}, \ \forall e\in E \label{pf_x_e}
\end{align}
\end{subequations}

The number of constraints in formulation \eqref{ef} grows exponentially with problem size, so we solve it via Benders Decomposition (BD).  At a typical iteration of BD, we consider the \emph{relaxed master problem} (RMP) \eqref{rmp}, which has the same objective as \eqref{ef} but involves only a small subset of the constraints in \eqref{ef}.  We briefly outline BD below.  For a detailed treatment of BD please refer to \cite{Benders1962}.

Let $j$ be the iteration counter and let the initial RMP be problem \eqref{rmp} without any constraints \eqref{rmp_f_cut}.  Let $x^j$ be a concatenation of the expansion variables $(x_e^j,x_g^j)$.

\begin{algorithm}[H]
\caption{\emph{Benders Decomposition} (BD)}
\begin{algorithmic}[1]
\State $j\gets 0$
\State solve RMP and let $x^j$ be the solution
\State \textbf{for} {$s=1,\cdots, S(k)$}
\State \hskip 0.4cm \textbf{if} {DSP$(x^{j},s)> \epsilon \sum_{i \in I} D_i$}
\State \hskip 0.8cm add feasibility cut $(\ref{dsp3})$ to RMP
\State \hskip 0.4cm \textbf{end if}
\State \textbf{end for}
\State \textbf{if} {$\forall \ s=1,\cdots,S(k)$, \ {DSP$(x^{k},s)\le \epsilon \sum_{i \in I} D_i$}}
\State \hskip 0.4cm $x^j$ is optimal (EXIT)
\State \textbf{else}
\State \hskip 0.4cm $j\gets j+1$ and GOTO step 2
\State \textbf{end if}
\end{algorithmic}
\end{algorithm}

By using a Benders reformulation, we are able to decompose the extremely large MINLP \eqref{ef} into a master problem and multiple subproblems (one for each contingency state).  In theory, this enables us to solve larger instances, which would not be possible by a direct solution of EF.  However, the extremely large number of contingency states makes direct application of Benders ineffective for large power systems and/or a non-trivial contingency budget (i.e., $k>1$). In the next section, we develop a custom cutting plane algorithm that evaluates all possible contingency states implicitly using a bilevel separation oracle.

\subsection{Cutting Plane Algorithm}
The size of most power systems in operations, with thousands of generating units and transmission elements, may preclude the direct solution of \eqref{ef}. Even using a decomposition algorithm (e.g. BD) may not be feasible because each contingency state must be considered explicitly. Our goal is to instead use a separation oracle that implicitly evaluates all contingency states and either identifies a violated one (a contingency with $k$ or less failures) or provides a certificate that no such contingency state exists. If such a violated contingency exists, we use this contingency to generate a Benders cut, as described in the previous section, for the RMP. If no such contingency exists, then the current capacity expansion $(x)$ is optimal and we terminate.

\subsubsection{Power System Inhibition Problem}
Given a capacity expansion decision $(x_e, x_g)$, the \emph{Power System Inhibition Problem} (PSIP) can be used  to determine the worst-case loss-of-load under any contingency with $k$ or fewer failures. In this bilevel program, the upper level decisions $(d_e, d_g)$ correspond to binary contingency selection decisions and the lower level decisions $(f,p, q,\theta)$ correspond to recourse power flow and load shedding decisions relative to this given contingency.

Note that, whereas in the prior model $\hat{d}^s$ was an input parameter, in this formulation, we are now selecting the elements of the contingency, with $d^s$ becoming a decision variable. For clarity of exposition, the superscript $s$ corresponding to variables $f, p, q$ and $\theta$ have been removed. PSIP is given as follows:
\begin{subequations}\label{psip}
\begin{align}
\hskip -0.1cm z(\hat x) = &\max_{d}  \  \min_{f,p,q,\theta} \quad   \sum_{i \in I}  q_i   \label{psip_obj}\\
\textrm{s.t.}  \quad & \sum_{e \in  E} d_e + \sum_{g \in G} d_g \leq k, \label{psip_budget}\\
	 (\alpha_i)  \quad & \sum_{g\in  G_i} p_g  + \sum_{e \in E_{.i}} f_{e} - \sum_{e \in E_{i.}} f_{e} +q_i = D_i,  \quad  \forall i \in  I  \label{psip_flow_consv}\\
	(\hat \beta_e) \quad &-B_{e}\big(\theta_{i_e} - \theta_{j_e}\big) + f_{e}  \leq M_e(1-\hat x_e +  d_e), \ \forall e \in  E  \label{psip_pa_ub}\\
	(\check \beta_e) \quad &B_{e}\big(\theta_{i_e} - \theta_{j_e}\big) - f_{e} \leq M_e(1-\hat x_e +  d_e), \ \forall e \in  E  \label{psip_pa_lb} \\
	(\delta_e) \quad &f_{e} \leq  F_{e} \hat x_e(1- d_e),  \quad \forall e \in  E    \label{psip_flow_ub}\\
	(\eta_e) \quad &-f_{e} \leq  F_{e} \hat x_e(1- d_e),  \quad \forall e \in  E    \label{psip_flow_lb}\\
	(\zeta_g) \quad & 0 \leq p_g \leq \overline{P}_g  \hat  x_g (1-d_g), \quad \forall g \in G  \label{psip_gen_ub}  \\
	(\lambda_i) \quad & 0 \leq q_i \leq D_i, \quad \forall i \in I \label{psip_lol_ub}
\end{align}
\end{subequations}

The objective \eqref{psip_obj} is to maximize the minimum loss-of-load. Given a contingency state defined by $(d_e, d_g)$, the objective of the power system operator (the inner minimization problem) is to determine the optimal power flow such that the loss-of-load is minimized. Constraint \eqref{psip_budget} is a budget constraint limiting the number of power system elements that can be in the contingency. Constraints \eqref{psip_flow_consv} are standard flow conservation constraints. Constraints \eqref{psip_pa_ub} and \eqref{psip_pa_lb} together enforce Kirchhoff's voltage law, if a transmission element is active. Constraints \eqref{psip_flow_ub}-\eqref{psip_flow_lb} are constraints associated with the capacity of each transmission element. Constraints \eqref{psip_gen_ub} limit the maximum capacity of each generating unit.  If a generating unit $g$ is NOT part of the contingency, that is $d_g = 0$, then the maximum capacity of the generating unit is enforced, if the unit was added ($x_g=1$).  Else, the power output of the generating unit must be zero.

The upper-level decisions of this bilevel program are to select a contingency, using the binary variables $(d_e, d_g)$, that maximizes the subsequent loss-of-load in the lower-level problem.

Bilevel programs like \eqref{psip} cannot be solved directly. One approach is to reformulate the bilevel program by dualizing the inner minimization problem. For fixed values of $d_e,d_g$, the inner minimization problem is a linear program that is always feasible. By using duality of linear programs, additional variables and disjunctive constraints, we can obtain equivalent an equivalent MILP formulation for \eqref{psip}, which we call the {\em Dual Power System Inhibition Problem} (D-PSIP). Bilevel programming approaches were used by \cite{SalWB09,Fan2011} to perform vulnerability analysis of power systems.

Next, we outline an algorithm for \emph{optimally} solving problem \eqref{ef} that combines a Benders decomposition with the aid of an oracle given by \eqref{psip}, which acts as a separation subroutine.  A given capacity expansion  $(x_e, x_g)$ is optimal if the oracle cannot find a contingency with $k$ or fewer failures that results in a loss-of-load above the allowable threshold $\epsilon \sum_{i \in I} D_i$. Whenever the oracle determines that the capacity expansion decision $(x_e,x_g)$ is not $N$-$k$-$\epsilon$ compliant,  it returns a contingency $(d_e, d_g)$ with $ \sum_{e \in E} d_e + \sum_{g \in G} d_g \leq k$.  This new contingency, given by $(d_e, d_g)$, results in a loss-of-load above the allowable threshold.

Let $j$ be the iteration counter, and let the initial RMP be problem \eqref{rmp} without any \eqref{rmp_f_cut} constraints and $x^j$ be the concatenation of  generation and transmission expansion variables.

\begin{algorithm}[H]
\caption{\emph{Cutting Plane Algorithm} (CPA)}
\begin{algorithmic}[1]
\State $j \gets 0$
\State solve RMP and let $x^j$ be the solution
\State solve D-PSIP$(x^{j})$ and let $d^j$ be the  solution
\State \textbf{if} { {D-PSIP$(x^j)> \epsilon \sum_{i\in I} D_i$} } 
\State \hskip 0.4cm solve DSP$(x^j,d^j)$
\State \hskip 0.4cm add feasibility cut \eqref{dsp3} to RMP \eqref{rmp}
\State \hskip 0.4cm $j \gets j+1$ and GOTO step 2
\State \textbf{else}
\State \hskip 0.4cm $x^j$ is optimal (EXIT)
\State \textbf{end if}
\end{algorithmic}
\end{algorithm}

At each iteration, either a contingency that results in loss-of-load above the allowable threshold is identified and a corresponding feasibility cut is generated and added to RMP, or the algorithm terminates with the current solution being optimal (if no failure contingency is found).

\section{Numerical Experiments}\label{sec4}
We implemented the above models and algorithms in C++ and CPLEX 12.1 via ILOG Concert Technology 2.9. All experiments were run on a machine with four quad-core 2.93G Xeon with 96G of memory.  For the following computational experiments, a single CPU and up to 8GB of RAM was allocated. The gap was set to be 0.1\% for CPLEX.

Our models and algorithms are tested on the IEEE-30-Bus and IEEE-57-Bus systems \cite{TestData}. For each power system, we consider five different contingency budgets $k$. Specifically, we limit the contingency cardinality $k$ to be at most zero, one, two, three, or four generating units or transmission elements in the power system. Altogether, we consider 10 instances.

Table \ref{tab1} allows us to compare the run times for the three different approaches. For each of the 10 instances, $m$ provides the number of distinct contingencies. Initially for each test system, we replicate a subset of existing generating units and transmission lines to create a \emph{starting} set of candidate elements. Using these candidate elements as a starting point, we iteratively solve the PSIP problem for $k=4$ and $\epsilon=0.05$ using CPA.  Using this process, we identify vulnerabilities in the power system and introduce additional candidate generation and transmission elements. We follow this method to create the augmented the IEEE-30-Bus and IEEE-57-Bus test systems for the computational experiments presented subsequently.

 \begin{table}[h]
\caption{Run times for different solution approaches}
\centering
\begin{tabular}{c c c c c c c}
\hline \hline
Test Systems 		&$m$ 		& $k$	&  $\epsilon$ 		& EF        		& BD        	& CPA           	\\[0.5ex]
\hline
IEEE-30-Bus		&	0		&	0	&	0 		    	& 0 			& 0 			&0			 \\
         				&	152		&	1	&	0 		    	& 76 			& 1 			& 2				    \\
            			&	$>11K$	&	2	&	0.05 			& x 			& 51 			& 15				\\
        				&	$>500K$	&	3	&	0.10 			& x 			& 2,382 		& 78				\\
        				&	$>21M$	&	4	&	0.20 			& x 			& x 			& 117				\\ [0.5ex]
\hline
IEEE-57-Bus     	&	0		&	0	&	0 		    	& 0 			& 0 			& 0				    \\
         				&	110		&	1	&	0 		    	& 27 			& 87 			& 74				    \\
            			&	$>5K$	&	2	&	0.05 			& x 			& 513 		& 15				\\
        				&	$>200K$	&	3	&	0.10 			& x 			& x 			& 23				\\
        				&	$>5M$	&	4	&	0.20 			& x 			& x 			& 24				\\ [0.5ex]
\hline
\end{tabular}
\label{tab1}
\end{table}

Table \ref{tab1} provides the run time (in CPU seconds) for each instance under the three different approaches. Note that the first approach, the extensive form (EF), can only solve the smallest of instances. This is because of the sheer size of the problem, in which, for each contingency, a full DCOPF problem must be embedded in the formulation. As the number of contingencies grows, this quickly becomes intractable.

The second approach, BD, bypasses this problem via a Benders decomposition, with corresponding delayed cut generation. However, this still suffers from the combinatorial growth in the number of contingency scenarios -- for each contingency, a subproblem (DSP) must be solved to check for violated feasibility cuts to add to the RMP. We see that larger problem instances can be solved, relative to EF, but the BD approach nonetheless cannot solve the largest problem instances.

In the CPA approach, we see that all instances of the problem can be solved, in all cases in under two minutes and frequently in only a few seconds. This is a result of the combination of the strength of the Benders cuts, enabling the problem to be solved in a very limited number of iterations, and also the fact that we are able to implicitly evaluate the contingencies in order to identify a violated contingency and then quickly find its corresponding feasibility cut by solving a single linear program (DSP).

Table \ref{tab2} provides us with further evidence to support this. For each instance, we see the total number of possible contingency scenarios $m$ and then the number of contingency scenarios for which corresponding feasibility cuts were actually generated (this is the total number of iterations itr). Clearly, it is a very tiny fraction of the possible number of contingencies, which is critical to the tractability of the approach. The remaining columns of this table breakdown the total run time by time spent on the three components of the algorithm -- the restricted master problem (RMP), which identifies a candidate network design $x$; the power system inhibition problem (PSIP), which identifies a contingency that cannot be overcome by the current network design; and the dual subproblems (DSP), which generates the feasibility cuts.

\begin{table}[h]
\caption{CPA runtime breakdown}
\centering
\begin{tabular}{c c c c c c c c}
\hline \hline
Test Systems 	&$m$ 		& $k$	&  $\epsilon$ 	&RMP        	&PSIP       	&DSP        	& itr		\\[0.5ex]
\hline
IEEE-30-Bus		&	0		&	0	&	0 		    & 0 		    & 0		        & 0				& 1		\\
         		&	152		&	1	&	0 		    & 0 			& 2 			& 0				& 6		\\
        		&	$>11K$	&	2	&	0.05 		& 0 			& 14 			& 0			    & 8		\\
        		&	$>500K$	&	3	&	0.10 		& 1 			& 78 			& 0			    & 14		\\
        		&	$>21M$	&	4	&	0.20 		& 1 			& 116		    & 0			    & 19		\\ [0.5ex]
\hline
IEEE-57-Bus    &	0		&	0	&	0 		    & 0 		    & 0		        & 0				& 1		\\
       		&	110		&	1	&	0 		    & 46			& 28 			& 0				& 103		\\
        		&	$>5K$	&	2	&	0.05 		& 7 			& 14 			& 0			    & 26		\\
        		&	$>200K$	&	3	&	0.10 		& 2 			& 21 			& 0			    & 39		\\
        		&	$>5M$	&	4	&	0.20 		& 1 			& 23		    & 0			    & 39		\\ [0.5ex]
\hline
\end{tabular}
\label{tab2}
\end{table}

\section{Conclusion}\label{sec5}
In this paper, we proposed models for TGEP with $N$-$k$--$\epsilon$ survivability constraints. Two algorithms are presented and tested on standard IEEE test systems. Computational results show the proposed custom cutting plane algorithm (CPA), using a bilevel separation program to implicit consider all exponential number of contingencies, significantly outperforms a standard Benders decomposition.  

The $k$ or fewer failures considered in this paper are assumed to happen simultaneously. In order to reflect practical operation situations, where failures may happen consecutively, new models that consider timing between system element failures will be needed. Additionally, unit commitment and de-commitment is not considered in our current model. Future research should address these considerations.

\vspace{5pt}
{
{\bf Acknowledgements.} This work was funded by the applied mathematics program at the United States Department of Energy and by the Laboratory Directed Research \& Development (LDRD) program at Sandia National Laboratories. Sandia National Laboratories is a multiprogram laboratory operated by Sandia Corporation, a wholly owned subsidiary of Lockheed Martin Corporation, for the United States Department of Energy's National Nuclear Security Administration under contract DE-AC04-94AL85000.
}

\newpage


\begin{thebibliography}{1}
\bibitem{NERC}
North American Electric Reliability Corporation (NERC), ``Standard TPL-001-1--system performance under normal conditions," February 2011. Available at http://www.nerc.com/docs/standards/sar/Project\_2006-02\_TPL-001-1.pdf

\bibitem{Pinar2010}
A. Pinar, J. Meza, V. Donde and B. Lesieutre, ``Optimization strategies for the vulnerability analysis of the electric power grid," {\em SIAM J. Optim.}, vol. 20, no. 4, pp.~1786--1810, 2010.

\bibitem{Bienstock2010}
D. Bienstock, A. Verma, ``The N-k problem in power grids: new models, formulations, and numerical experiments," {\em SIAM J. Optim.}, vol. 20, no. 5, pp.~2352--2380, 2010.

\bibitem{Arroyo2010}
J.M. Arroyo, ``Bilevel programming applied to power system vulnerability analysis under multiple contingencies," {\em IET Gener Transm. Distrib.}, vol. 4, no. 2, pp. 178--190, 2010.

\bibitem{SalWB04}
J. Salmeron, K. Wood, R. Baldick, ``Analysis of electric grid security under terrorist threat," {\em IEEE Trans. Power Syst.}, vol. 19, no. 2, pp. 905--912, 2004.

\bibitem{SalWB09}
J. Salmeron,  K. Wood, R. Baldick, ``Worst-case interdiction analysis of large-scale electric power grids,"  {\em IEEE Trans. Power Syst.}, vol. 24, no. 1, pp. 96--104, 2009.


\bibitem{Fan2011}
N. Fan, H. Xu, F. Pan, P.M. Pardalos, ``Economic analysis of the N-k power grid contingency selection and evaluation by graph algorithms and interdiction methods," {\em Energy Syst.}, vol. 2, no. 3--4, pp. 313--324, 2011.

\bibitem{Street2011}
A. Street, F. Oliveira, J.M. Arroyo, ``Contingency-constrained unit commitment with $n$-$K$ security criterion: A robust optimization approach," {\em IEEE Trans. Power Syst.}, vol. 26, no. 3, pp. 1581--1590, 2011.



\bibitem{Romero2011}
N. Romero, N. Xu, L.K. Nozick, I. Dobson, D. Jones, ``Investment planning for electric power systems under terrorist threat," {\em IEEE Trans. Power Syst.}, DOI: 10.1109/TPWRS.2011.2159138, 2011.


\bibitem{Carrion2007}
M. Carrion, J.M. Arroyo, N. Alguacil, ``Vulnerability-constrained transmission expansion planning: A stochastic programming approach," {\em IEEE Trans. Power Syst.}, vol. 22, no. 4, pp. 1436--1445, 2007.

\bibitem{Choi2007}
J. Choi, T.D. Mount, R.J. Thomas, ``Transmission expansion planning using contingency criteria," {\em IEEE Trans. Power Syst.}, vol. 22, no. 4, pp. 2249--2261, 2007.

\bibitem{Moulin2010}
L. Moulin, M. Poss, C. Sagastiz\'{a}bal, ``Transmission expansion planning with re-design," {\em Energy Syst.}, vol. 1, no. 2, pp. 113--139, 2010.

\bibitem{Bent2011}
R. Bent, A. Berscheid, A.L. Toole, ``Generation and transmission expansion planning for renewable energy integration," in {\em Proc. of Power Syst. Comp. Conf. (PSCC)}, Stockholm, Sweden, August 2011.


\bibitem{Jin2011}
S. Jin, S.M. Ryan, J.-P. Watson, D.L. Woodruff, ``Modeling and solving a large-scale generation expansion planning problem under uncertainty," {\em Energy Syst.}, vol. 2, no. 3--4, pp. 209--242, 2011.



\bibitem{Benders1962}
J.F. Benders, ``Partitioning procedures for solving mixed-variables programming problems," {\em Numerische Mathematik}, vol. 10, pp.~237--260, 1962.

\bibitem{TestData}
IEEE reliability test data, [online] Available at http://www.ee.washing-\\
ton.edu/research/pstca/.


\end{thebibliography}
\end{document}